\documentclass[12pt]{amsart}
\newtheorem{theorem}{Th\'eor\`eme}
\newtheorem{proposition}{Proposition}
\newtheorem{definition}{Definition}
\newtheorem{notation}{Notation}
 \newtheorem{corollary}{Corollaire}
\newtheorem{lemma}{Lemme}

\input{amssym.def}
\begin{document}
\title{applications  harmoniques et hyperbolicit\'e
de domaines tubes }
\author{Jean-Jacques Loeb}
\maketitle

\begin{abstract}
An application of the Zalcman renormalization
theorem   to  harmonic
functions shows that the  limit functions are  nonconstant affine. 
Extensions of this method are given for maps 
with values in a torus or in a complex Lie groups. 
As an application, we give  criteria of Kobayashi 
hyperbolicity for tubes in ${\Bbb C}^2$.  

.\\

{\bf Mathematics subject classification}(2000)
 32A07, 31CO5\\

{\bf Keywords} tube domain. Renormalization. Harmonic maps. 
.    
\end{abstract}

\section{\bf Introduction}
Dans leur article \cite{b}, F. Berteloot et 
J. Duval d\'emontrent le th\'eor\`eme suivant: 

\begin{theorem} Etant donn\'e une fonction holomorphe
enti\`ere $f$  non constante qui ne s'annule 
pas, il existe une suite $A_n$ 
de nombres positifs et
une suite  $B_n$ de nombres complexes 
telles que la suite des fonctions 
$f(A_nz +B_n)$ tende uniform\'ement
sur tout compact vers une fonction de 
la forme $Ce ^{dz}$ avec $C$ et $d$ non nuls. 
\end{theorem} 

Ils prouvent ce th\'eor\`eme en utilisant le 
th\'eor\`eme de renormalisation de 
Zalcman \cite{z}. Un th\'eor\`eme de 
renormalisation similaire
\`a celui de Zalcman  avait \'et\'e
\'egalement obtenu par Brody  \cite{br}
qui utilisait des homographies au lieu de transformations affines. 
Comme l'avaient remarqu\'e Berteloot et Duval, 
le petit th\'eor\`eme de Picard est 
une cons\'equence imm\'ediate de leur r\'esultat  
moyennant l'utilisation d'un lemme   classique d'Hurwitz.
Dans  leur article, les auteurs donnaient 
une g\'en\'eralisation de leur th\'eor\`eme 
\`a la dimension sup\`erieure, ce qui leur permettait 
d'\'etablir l'hyperbolicit\'e du compl\'ementaire 
de certaines courbes dans le plan projectif complexe. 
Ils retrouvaient aussi avec une m\'ethode 
\'el\'ementaire la g\'en\'eralisation du th\'eor\`eme 
de Picard par Mark Green, 
\`a savoir que le compl\'ementaire de $2n+1$ 
hyperplans
en position g\'en\'erale dans ${\Bbb P}^n$  est hyperbolique (ici ceci signifie 
essentiellement qu'une application holomorphe  enti\`ere qui prend ses 
 valeurs dans ce compl\'ementaire est constante). 

Comme point de d\'epart et motivation de 
n\^otre \'etude, notons le fait suivant:

\begin{proposition} 
On a l'\'equivalence entre le th\'eor\`eme 
pr\'ec\'edent et le r\'esultat suivant: 

Pour toute fonction harmonique $g$ 
d\'efinie sur ${\Bbb C}$ non constante, 
il existe une suite de positifs $A_n$ 
et une suite de complexes $B_n$ tels que 
la suite $g(A_nz+B_n)$ tende uniform\'ement 
sur tout compact vers une fonction affine non constante. 
\end{proposition} 

La preuve de ce r\'esultat est comme suit.
Tout d'abord si $g$ est harmonique sur ${\Bbb C}$, 
on l'\'ecrit comme la partie r\'eelle 
d'une fonction enti\`ere $h$ et de ce 
fait $g=\ln \vert f\vert$, o\`u $f=e ^h$. 
On d\'emontre ainsi que le th\'eor\`eme 
de Berteloot et Duval implique 
le r\'esultat analogue sur les fonctions harmoniques. 
Dans l'autre direction, si $f$ est une fonction 
enti\`ere ne 
s'annulant pas, on commence par appliquer le r\'esultat 
de renormalisation \`a  la fonction harmonique 
$g=\ln \vert f\vert$. On a donc une suite $g(A_nz+B_n)$ 
qui tend vers une fonction affine non constante. On utilise 
la conjugu\'ee $l$ de $g$ et les relations de Cauchy-Riemann 
pour montrer que pour une d\'etermination  $\ln f$ 
  de $f$, il existe une suite de r\'eels $c_n$ 
tels que $\ln f(A_nz+B_n)+ic_n$ ait une limite
qui est n\'ecessairement une fonction holomorphe affine non constante.  
On peut maintenant s'arranger quitte \`a extraire, 
pour que la suite $\ln f(A_nz+B_n)$ ait une limite 
modulo $2i\pi {\Bbb Z}$. Ceci permet de conclure. \\

Notons qu'un lemme 
de renormalisation est \'etabli en dimension deux dans \cite{s}
pour les fonctions harmoniques sans qu'il soit donn\'e  
de r\'esultat pr\'ecis sur la limite. 

N\^otre article s'organise comme suit: 

Dans la premi\`ere partie, nous  g\'en\'eralisons le 
r\'esultat concernant les fonctions harmoniques 
\`a toutes les dimensions. Nous suivons la d\'emarche 
de \cite{b} et  
nous utilisons les in\'egalit\'es classiques de Harnack 
pour travailler sur les fonctions harmoniques. 
Nous donnons ensuite des crit\`eres de 
normalit\'e pour des familles de fonctions harmoniques 
ainsi qu'une caract\'erisation intrins\`eque des fonctions affines 
parmi les fonctions harmoniques. 

Dans la seconde partie, nous faisons 
le lien entre applications harmoniques 
et tubes dans ${\Bbb C}^n$. Des crit\`eres d'hyperbolicit\'e 
(au sens de Kobayashi) sont donn\'es. On montre ainsi qu'il existe
des tubes hyperboliques et dont l'enveloppe d'holomorphie est 
tout ${\Bbb C}^2$. Il nous a paru int\'eressant 
de faire quelques remarques en compl\'ement 
  sur les images d'applications enti\`eres  harmoniques.

Dans la derni\`ere partie nous montrons que les 
fonctions (ou applications) solutions d'op\'erateurs elliptiques 
admettent une renormalisation \` a constante 
pr\`es en des fonctions (ou applications) affines non constantes. 
Une g\'en\'eralisation naturelle est donn\'e 
pour les groupes de Lie complexes. 

Pour les propri\'et\'es \'el\'ementaires sur les 
fonctions harmoniques utilis\'ees ici, nous renvoyons le 
lecteur \`a \cite{h}.   

\section{\bf Renormalisation des fonctions harmoniques}

\begin{notation} Dans cette partie, on d\'esigne par $U$ 
un domaine dans ${\Bbb R}^m$ et on note par  $H(U)$ l'espace des 
fonctions harmoniques sur $U$. Une fonction de 
$H({\Bbb R}^m)$ sera dite harmonique enti\`ere.  

Sur $[-\infty, +\infty]$ on met la distance $d$ 
d\'efinie par: $d(x,y)=\vert \ S( x) -  S( y) \vert $, 
avec $S(x):= \arctan (\sinh  x)$.  Cette 
distance d\'efinit la topologie usuelle. Elle est inspir\'ee de 
la m\'etrique de Fubini-Study.  

Pour une fonction $f$ diff\'erentiable 
sur $U$, on pose $\tilde{f} (x):= \vert (S\circ f)'(x)\vert =
s\circ f(x)\vert f'(x)\vert $ o\` u on a not\'e $s$ 
la d\'eriv\'ee de $S$ et $f'$ le gradient de $f$.
Explicitement: $\tilde{f}
 (x)= \frac{\vert f'(x)\vert }{\cosh f(x)}$ 
\end{notation}

Dans la suite, nous utiliserons un 
lemme de renormalisation de base donn\'ee par plusieurs auteurs.
La version que nous donnons,  
d\^ue \`a F. Berteloot, 
est similaire \`a  de celle de 
Gromov (voir \cite{g}). 

\begin{lemma} 
Soit $(V,d)$ un espace m\'etrique 
complet et $\phi$ une fonction 
sur $V$ \`a valeurs positives, localement born\'ee. 
On fixe $\tau >1 $, $\epsilon >0$ et $p\in V$ tel que $\phi (p)>0$. 
Alors il existe $q\in V$ tel que: 

1. $d(p,q)\leq \frac{\tau}{\epsilon \phi (p)(\tau -1)}$.  

2. $\phi (q)\geq \phi (p)$

3. $\phi (x)\leq \tau \phi (q)$ si $d(x,q)\leq \frac{1}{\epsilon 
  \phi (q)}$.\end{lemma}

La preuve de ce lemme est  \'el\'ementaire. 
Elle consiste \`a construre une suite 
$p_n$ en partant de $p_0=p$. Si la condition 3. 
est satisfaite en rempla\c{c}ant $q$ 
par $p_0$, on s'arr\`ete et on prend $q:=p_0$. 
Sinon on choisit  $p_1$ dans la boule ferm\'ee
de centre $p_0$ et de rayon 
$\frac{1}{\epsilon \phi (p_0)}$ 
avec $\phi (p_1)>\tau \phi (p_0)$. On raisonne alors sur $p_1$ comme sur $p_0$. 
Il est alors facile de voir en utilisant 
l'hypoth\`ese $\phi$ localement born\'ee que pour un certain $p_n$ 
la condition 3. est satisfaite, ainsi que les deux autres. 

Le lemme suivant, cons\'equence du lemme pr\'ec\'edent 
est ce qui nous int\'eressera par la suite.

\begin{lemma} 
 Soit $f_n$ une suite de fonctions 
diff\'erentiables sur $U$. On suppose que pour un certain 
point $r$ de $U$, il existe une suite 
$r_n$  qui tend vers $r$ et telle que  $\tilde{f_n} (r_n)$ tend vers 
l'infini. Alors il existe deux suites $a_n>0$ et $b_n\in {\Bbb R}^m$
tendant vers $0$ et $r$ respectivement  et telles que si on pose: 
$g_n(z):= f_n(a_n z +b_n)$, on a:  d'une part
 $\tilde{g_n}
 (0)=1$ et d'autre part   
 pour tout compact $K$ de ${\Bbb R}^m$, il existe 
une suite de nombres $\epsilon _n >0$ tendant vers $0$ 
et telle que: $\tilde{g_n} \leq 1+\epsilon _n$ sur $K$.\end{lemma}

\emph{preuve} On particularise d'abord le lemme pr\'ec\'edent 
en prenant 
$\epsilon:=\phi (p)^{-\frac{1}{3}}$, 
$\tau:=1+\epsilon$.  La condition 1. devient: $d(p,q)\leq (1+\phi (p)^
{-\frac{1}{3}})\phi(p)^{-\frac{1}{3}}$.  
  On pose ensuite:  $p:=r_n$ et $\phi:=\tilde {f_n}$
et on choisit pour $V$ du lemme pr\'ec\'edent une boule ferm\'ee non
triviale de centre $r$ contenue dans $U$. On applique \`a cette situation 
le lemme particularis\'e. Pour chaque $n$, le lemme nous donne 
un \'el\'ement $q$  not\'e $q_n$. En utilisant la 
condition 1. et l'hypoth\`ese $\tilde {f_n}(r_n)$
tend vers l'infini, on voit que $q_n$ tend vers $r$. 
Pour conclure, on prendra pour suite $b_n$ la suite des $q_n$ et 
pour $a_n$ la suite des $\frac{1}{\tilde {f_n} (q_n)}$. En utilisant
toujours l'hypoth\`ese, la condition 3. nous montre que 
la suite des fonctions $g_n$ construite \`a partir des $a_n$ et $b_n$ 
satisfait bien aux conditions requises. \\

\emph{D\'efinition}  On dira qu'une telle suite $g_n$ 
est une renormalisation de  la suite $f_n$ en $r$. \\ 

Cons\'equence du lemme: En utilisant l'in\'egalit\'e des 
accroissements finis, on voit que la famille des 
fonctions $g_n$ est une famille \'equicontinue \`a valeurs 
dans ${\Bbb R}$. Ceci implique via le th\'eor\` eme 
d'Ascoli que si on voit les $g_n$ \`a valeurs dans 
le compact $[-\infty, +\infty]$ on peut en extraire 
une sous-suite qui converge uniform\'ement sur tout 
compact vers $g$ \`a valeurs dans $[-\infty,+\infty]$.

{\bf Remarque}  Le principe de renormalisation 
sera utilis\'e plusieurs fois dans la suite 
sans que nous donnions toutes les \'etapes, 
comme nous les avons donn\'ees ici.

Enon\c{c}ons  le th\'eor\`eme principal: 

\begin{theorem} Soit une suite $f_n$ de fonctions
de $H(U)$. On suppose que pour un point $p$, la suite des 
$\tilde{f_n}
 (p)$ tende  vers l'infini. Alors on peut renormaliser en $p$
la  
suite des $f_n$ en une suite $g_n$ 
telle qu'une suite extraite  tende  uniform\'ement sur tout
compact de ${\Bbb R}^m$ vers une fonction affine non constante. 
\end{theorem} 

\begin{corollary} Soit $f$ une fonction harmonique
enti\`ere  non constante.  
Alors il existe une suite de positifs $A_n$ et une suite
d'\'el\'ements $B_n$ de ${\Bbb R}^m$ telle que $f(A_nz +B_n)$ 
tende uniform\'ement sur tout compact vers une fonction affine non
constante. \end{corollary}

\emph{preuve du corollaire}  Comme $f$ est non constante, 
on peut trouver un point $p$ tel que $f'(p)$ soit non nul. 
On applique alors le th\'eor\`eme \`a la suite 
$f_n(z):= f(p+ nz)$. \\

{\bf Remarque} Dans le cas o\`u $f$ est la
 est la  partie 
r\'eelle d'un polyn\^ome holomorphe, il est 
facile de voir directement l'existence des $A_n$ 
et des $B_n$. On peut aussi v\'erifier 
dans le cas g\'en\'eral  qu'on peut choisir la suite des $A_n$ 
born\'ee. \\

\emph{preuve du th\'eor\`eme} La preuve du th\'eor\`eme va 
se faire  en deux temps: 
On prouve d'abord que la suite renormalis\'ee $g_n$ 
tend apr\`es extraction vers $g$ \`a valeurs finies. 
L'ellipticit\'e du Laplacien permet d'en d\'eduire imm\'ediatement 
que $g$ est harmonique avec $\tilde{g}
 \leq 1$ et $\tilde{g} (0)=1$.  
Cette derni\`ere \'egalit\'e montre que $g$ est non constante. 
On exploite ensuite la condition $\tilde{g} \leq 1$ 
pour conclure.

La finitude de $g$  se d\'eduit 
imm\'ediatement  du lemme suivant: 

\begin{lemma}  Soit $g_n$ une suite de fonctions de 
$H(U)$ qui tend
uniform\'ement sur tout compact  vers $g$ \`a valeurs dans
$[-\infty,+\infty]$. Alors si $g$ vaut $+\infty$ ou $-\infty$ 
en un point, elle va \^etre \'egale \`a cette 
valeur partout. De plus $\tilde{g}$   est alors identiquement
nulle.\end{lemma} 

\emph{preuve} On utilise les in\'egalit\'es de Harnack
(cons\'equences de la formule int\'egrale de Poisson)
dont  voici un \'enonc\'e:  Etant donn\'e une boule ouverte de rayon $2R$ 
et de centre $p$, il existe une  constante $A$  
strictement positive telle que pour toute fonction harmonique 
$f$ positive dans cette boule et pour tout $x$ et $y$ 
dans la boule de centre $p$ et de rayon $R$, on ait: 
$Af(x)\leq f(y)$ (voir \cite{h}).

Revenons au lemme  et supposons par ex $g(p)=+\infty$. 
Alors par convergence uniforme les $g_n$ seront 
positives dans un voisinage de $p$ pour $n$ assez grand. 
On va donc avoir pour $x$ dans un  voisinage de $p$:  $g_n(x)\geq A
g_n(p)$. Ceci implique que $g_n$ tend vers $+\infty$ 
dans un voisinage de $p$. Par convergence uniforme, l'ensemble des 
points $x$ o\`u $g_n(x)$ tend vers $+\infty$ est ferm\'e. 
La connexit\'e de $U$ permet de conclure que $g_n$ tend partout vers
$+\infty$. Pla\c{c}ons nous alors dans cette situation et fixons $p$ 
dans $U$. Les in\'egalit\'es de Harnack nous montrent que 
la suite des $v_n(x):= (1/c_n)  g_n(x)$
avec $c_n=g_n(p)$  est une suite born\'ee dans 
un voisinage de $p$. Par normalit\'e,
 on peut donc en extraire une sous suite qui
converge dans un voisinage vers une fonction qui sera  strictement 
positive toujours d'apr\`es Harnack. Pour montrer que 
les $\tilde{g_n}
 (p)$ tendent  vers z\'ero, on \'ecrit: 
$\tilde{g_n}(p)=  c_n\vert v_n'(p)\vert s(c_nv_n(p))$ 
et on utilise la propri\'et\'e que $xs(x)$ tend vers 
$0$ quand $s(x)$ tend vers l'infini. \\

Passons \`a la seconde partie de la preuve. Par passage \`a la limite 
et ellipticit\'e, on a maintenant une fonction limite $g$ harmonique 
telle que $\tilde{g} \leq 1$. Cette partie est similaire 
\`a la preuve donn\'ee dans \cite{b}. On remarque 
que pour $g$ harmonique, le laplacien de $\ln (\cosh g)$ 
vaut $\frac{\vert g'(x)\vert ^2}{(\cosh g(x) )^2}$ 
soit $\tilde{g} ^2$. On en d\'eduit que $c\vert x\vert ^2 -\ln (\cosh
g)$ est sous-harmonique en prenant $c=1/{2m}$. 
De ce fait si on note $M(F,r)$ l'int\'egrale d'une fonction $F$ 
sur la sph\`ere de centre $0$ par rapport \`a la mesure $\mu$ 
standard de masse un  sur la  sph\`ere, on a d'abord
$M(\ln \cosh g,r)\leq cr^2$ puis en tenant compte de la 
relation $\vert x\vert \leq  (\ln \cosh x)+d $ avec $d=\ln 2$ 
 on en d\'eduit: 
$M(\vert g\vert ,r)\leq cr^2+d$. Il est alors classique que ceci 
implique que $g$ est un polynome de degr\'e au plus deux. 
La m\'ethode consiste \`a d\'evelopper $g$ en s\'erie de 
polyn\^omes harmoniques homog\`enes
$P_k$. Les relations d'orthogonalit\'e 
des $P_k$ permettent d'en d\'eduire: 
$M(gP_k,r)= M(P_k^2,r)$. On fait alors tendre $r$ 
vers l'infini  en tenant compte de   l'homog\'en\'eit\'e 
des $P_k$ et de l'in\'egalit\'e sur les $M(g,r)$
pour conclure. La nullit\'e de la partie quadratique de $g$  
r\'esulte de 
la normalit\'e de la famille des 
$g_t(z):= g(z+t)$ pour $t\in {\Bbb R}^m$. (Pour la normalit\'e, voir la
section suivante et la proposition 2). 
En effet on \'ecrit: $g(z)=Q(z)+ L(z)$ avec $Q$ partie quadratique 
homog\`ene et $L$ partie affine. Si on suppose  $Q$  harmonique
non nulle, il existe
un 
vecteur $v$ isotrope qui n'est pas  dans le noyau de la 
forme bilin\'eaire $B$ associ\'ee. On consid\'ere alors 
la suite des $g(z+nv)= n(B(z,v)+ L(v))+L(z)$ 
et comme $B(z,v)+ L(v)$ peut prendre toutes les valeurs r\'eelles, 
il n'existe pas de sous suite extraite qui converge. 
Ceci termine la preuve du th\'eor\`eme.\\

\noindent Esquissons  une preuve directe 
du corollaire 
analogue \`a la premi\`ere preuve de \cite{b} . 
 Soit  $F$ harmonique enti\`ere non constante. 
Supposons par exemple $F'(0)$ non nulle. 
On choisit alors une suite $R_n$ de telle sorte 
que pour la suite $F_n(z)= F(R_nz)/n$ on ait: 
$\lim\tilde{F_n}
 (0)=\infty$. On renormalise alors cette 
suite en une suite $g_n$ comme dans la premi\`ere partie de la 
preuve du th\'eor\`eme. La limite $g$ de la suite extraite des 
$g_n$ est harmonique enti\`ere non constante . 
Il est classique qu'elle  s'annule en un point $p$ (\cite{h}) et   
comme elle est harmonique, on peut aussi 
supposer qu'en ce point, le gradient ne s'annule pas (voir lemme suivant).
Supposons par exemple  $p=0$. On conclut alors en consid\'erant la suite 
des $ng_n(z/n)$ et en faisant tendre $n$ vers l'infini.\\

Montrons le lemme dont on a eu  besoin. 

\begin{lemma} Soit $f\in H(U)$ non identiquement nulle. 
Alors si $f$ s'annule en un point $p$, il existe un point $q$ 
o\`u $f$ s'annule et tel que $f'(q)$ ne s'annule pas. \end{lemma}

\emph{preuve} D'apr\'es le principe du maximum, $f$ 
prend des valeurs strictement positives et strictement
n\'egatives. Donc l'ensemble des z\'eros de $f$ qui est un ensemble 
analytique r\'eel,  est une hypersurface. On se place alors en un point lisse
de cette  hypersurface et on applique un  lemme
classique  de Hopf pour conclure
\cite {h}.

\subsection{Crit\`eres de  normalit\'e}

Ces remarques sont fortement  inspir\'ees 
par l'article de Zalcman 
 \cite{z} qui tra\^ite le cas m\'eromorphe \`a une variable. 

\begin{definition}
(voir \cite{s} pour la dimension deux)

 On dira qu'une famille $E$ 
de fonctions de $H(U)$ est normale 
si de toute suite de $E$, on peut extraire une sous suite 
convergeant uniform\'ement sur tout compact soit 
vers une fonction de $H(U)$ soit vers $+\infty$ 
soit vers $-\infty$. \end{definition} 

On a la proposition suivante qui 
donne un \'equivalent du th\'eor\`eme de 
Marty pour le cas harmonique ( voir \cite{z}).

\begin{proposition} Une famille E est normale 
si et seulement si pour tout compact $K$, il existe 
une constante positive $M_K$ tel que 
pour tout $f\in E$, on ait: $\tilde{f}
 (x)  \leq M_K$
pour $x\in K$. \end{proposition} 

\emph{preuve} Si on a la condition de la proposition  sur les 
d\'eriv\'ees au sens du tilde, alors la famille est \'equicontinue 
et le th\'eor\`eme d'Ascoli joint au lemme 3.  
 montre la normalit\'e de la famille. 
Dans l'autre direction, supposons que la condition sur
les d\'eriv\'ees ne soit pas remplie. Il existe alors un 
compact $K$ et une suite $f_n$ de $E$ 
telle que le $\sup$ des  $\tilde
{f_n} $ tende vers l'infini sur $K$. 
Si la suite \'etait normale, quitte \`a extraire, 
on pourrait  supposer que les $f_n$ tende uniform\'ement sur tout 
compact de $U$ vers 
$f\in H(U)$ ou vers $+\infty$ ou vers $-\infty$.
Ceci est contradictoire avec  
ce qui a \'et\'e dit pr\'ecedemment sur 
les relations entre les limites de $f_n$ 
et de $\tilde{f_n}$. 

On a le crit\`ere g\'en\'eral de normalit\'e suivant: 

\begin{theorem} Une famille $E$ dans $H(U)$ 
est normale si elle satisfait \`a l'un des deux crit\`eres suivants: 

1. On fixe  $a\in {\Bbb R}$.
On consid\`ere la famille 
$E$ des $f\in H(U)$ v\'erifiant:   pour tout compact $K$, 
il existe $M_K>0$ tel que $\vert f'(x)\vert$ soit major\'ee 
par $M_K$ sur $f^{-1}(a)\cap K$ .  

2. On fixe une fonction $l$ sur ${\Bbb R}$ \`a valeurs dans 
$[0,\infty]$ finie en au moins un point. 
On consid\`ere la famille $E$ des $f\in H(U)$ 
v\'erifiant:  $\vert f'\vert \leq l\circ f$ \end{theorem}

\emph{preuve} 1. On suppose la famille non normale. 
Il existe d'apr\`es la 
proposition pr\'ec\'edente,   une suite $f_n$ dans 
$H(U)$ et une suite $p_n$ qui tend vers $p$ dans 
$U$ tel que $\tilde{f_n}
 (p_n)$ tende vers l'infini. 
Le lemme de renormalisation  
permet de renormaliser $f_n$ en $g_n$ 
au point $p$. Pour la suite, nous avons simplement 
besoin de savoir que la suite $g_n$ (modulo extraction) 
converge vers $g$ harmonique  enti\`ere  et 
non constante. En utilisant la relation: 
$g_n'(z)=a_n f_n'(a_nz+b_n)$ , pour $z\in K\cap f^{-1}(a)$ on a:  
$\vert g_n'(z)\vert \leq a_nM_K$ pour $n$ assez grand. 
On en d\'eduit que $g'$ s'annule sur $g^{-1}(a)$, ce qui 
est contradictoire avec le lemme pr\'ecedent. 

2. On choisit $b$ tel que $l(b)$ soit fini et on applique 
le 1.  avec les $M_K=l(b)$.  

\subsection{Fonctions de Brody}
 
En s'inspirant du cas m\'eromorphe, on dira que 
 $f$ harmonique enti\`ere  est une fonction de Brody si elle
v\'erifie $\tilde{f}
 \leq M$
pour un certain $M>0$.  La preuve du th\'eor\`eme 
pr\'ec\'edent montre qu'une fonction harmonique de Brody est en fait 
une fonction affine. Nous allons donner un th\'eor\`eme 
qui donne une caract\'erisation plus intrins\`eque 
des fonctions de Brody et qui a un analogue dans 
le cas m\'eromorphe. 

\begin{theorem} 
1. Une fonction harmonique enti\`ere  $f$  
est de Brody ( i.e affine) si et seulement si 
la famille des fonctions $f_t(x):=f(x+t)$ index\'ee 
par $t\in {\Bbb R}^m$ forme  une famille normale. 

2. Soit $l$ une fonction de ${\Bbb R}$ dans  
$[0,+\infty]$ finie en au moins un point. 
Alors les fonctions harmoniques $f$ 
sur ${\Bbb R}^m$ v\'erifiant $\vert f'\vert \leq l(f)$ 
sont de Brody (i.e. affines). \end{theorem} 

\emph{preuve} 
1. Supposons la famille $f_t$ normale. Alors 
d'apr\`es la proposition 2., il existe $M>0$ tel que
pour tout $t$, on a: $\tilde{f_t}
 (0)\leq M$ . 
Ce qui signifie $\tilde{f} \leq M $ et donc 
$f$ est affine. 
La r\'eciproque est claire dans n\^otre situation o\`u on sait 
que les fonctions de Brody sont affines. Mais 
on peut aussi donner un argument direct  
en remarquant simplement que si $\tilde{f}\leq M$, 
alors il en est de m\^eme de ses  translat\'ees 
qui forment alors 
une sous-famille de  la famille normale  des $g\in H({\Bbb R}^m)$ 
tels que $\tilde{g}\leq M$  

2. La famille $E$ des $f$ harmoniques enti\`eres telles que: 
$\vert f'\vert \leq l(f)$ est une famille normale d'apr\`es 
le th\'eor\`eme pr\'ec\'edent. On remarque ensuite que si $f\in E$, il 
en est de m\^eme des $f_t$. La famille des $f_t$ 
est normale et donc $f$ est affine d'apr\`es le 1.

\section{\bf Lien avec la g\'eom\'etrie des  tubes}

\emph{d\'efinition} On dira que $\Omega$ est un  tube (ou domaine tube)
de ${\Bbb C}^n$ s'il s'\'ecrit sous la forme $\omega +i{\Bbb R}^n$
avec $\omega$ domaine de ${\Bbb R}^n$. 

Citons \`a propos des  tubes un th\'eor\`eme classique 
de S. Bochner. \cite{h}

\begin{theorem} L'enveloppe d'holomorphie 
d'un domaine tube con\"{i}ncide avec son enveloppe convexe. 
\end{theorem} 

De ce th\'eor\`eme, on d\'eduit ais\'ement: 

\begin{corollary} Pour un domaine tube $\Omega$  
pseudo-convexe, on a \'equivalence entre:

1. $\Omega$ est biholomorphe \`a un domaine born\'e. 

2. $\Omega$ est hyperbolique (au sens de Kobayashi).  

3. $\omega$ ne contient pas de droite 
affine. 
\end{corollary}

\emph{preuve du corollaire} Par le th\'eor\`eme de 
Bochner, $\Omega$ est convexe. Le 3. signifie 
que $\Omega$ contient une droite affine complexe. 
L'\'equivalence de 1. et 3. est bien connue 
dans ce cas  (\cite{ba}). Elle implique 
\'evidemment l'\'equivalence commune avec 2. 

On utilise ce r\'esultat pour montrer: 

\begin{proposition} Un  tube $\Omega$ 
est biholomorphe \`a un domaine born\'e 
si et seulement si l'enveloppe convexe de $\omega$ 
ne contient pas de  droite affine. \end{proposition}

\emph{preuve} La condition suffisante  est 
\'evidente d'apr\`es le corollaire pr\'ec\'edent. 
Montrons que la condition est n\'ecessaire 
et supposons  que l'enveloppe convexe  $\omega '$ 
de $\omega$ 
contienne une droite $L$ . Alors par chaque point 
de $\omega '$ il passe une droite parall\`ele 
\`a $L$ et contenue dans $\omega'$. On peut alors choisir 
$x$ et $y$ distincts dans  $\omega$ et contenus 
dans une droite $D$ de $\omega '$. Ceci va impliquer 
que pour toute fonction holomorphe born\'ee $f$ 
sur $\omega$, on a: $f(x)=f(y)$, et donc $\Omega$ 
ne peut \^etre biholomorphe \`a un domaine born\'e.
L'assertion sur $f$ s'obtient comme suit: 
On peut supposer $f$ \`a valeurs dans le disque 
unit\'e. Cette fonction va se prolonger en une fonction holomorphe
$F$ 
de  l'enveloppe d'holomorphie (qui est aussi l'enveloppe convexe)
toujours \`a valeurs dans le disque unit\'e. Par Liouville, 
cette fonction va \^etre constante sur la droite complexe $D+iD'$, 
o\`u $D'$ est la droite vectorielle parall\`ele \`a $D$. En
particulier $f(x)=f(y)$. \\

\noindent {\bf remarque}  La preuve pr\'ec\'edente montre en fait que 
pour un domaine tube, on a l'\'equivalence entre: 
Le domaine est biholomorphe \`a un domaine born\'e, 
la pseudo-m\'etrique de Carath\'eodory 
est une m\'etrique, l'enveloppe convexe de $\omega$ 
ne contient pas de droite affine. \\

N\^otre int\`eret  concernant  l'hyperbolicit\'e 
va donc se porter sur les domaines tubes  non pseudo-convexes
dont l'enveloppe convexe contient des droites affines.
On introduit les d\'efinitions suivantes: 

\noindent \emph{d\'efinitions} 

\noindent 1. Une application harmonique (voir par ex. \cite{ah}) 
d'un domaine $U$ de ${\Bbb R}^n$ 
\`a valeurs dans ${\Bbb R}^m$ est 
une application dont les fonctions coordonn\'ees 
sont harmoniques.

\noindent 2. Dans le cas o\`u $U={\Bbb R}^n$, 
on parlera d'application harmonique enti\`ere. 

\noindent 3. Un domaine $V$ de 
${\Bbb R}^m$ est dit $n$- Brody-hyperbolique
(par rapport aux applications harmoniques)   
s'il n'existe pas d'application 
harmonique enti\`ere non constante 
de ${\Bbb R}^n$ \`a valeurs dans $V$.
Dans la suite, nous serons surtout int\'eress\'es 
par le cas $n=2$.  

\noindent Un domaine  $V$ de ${\Bbb R}^m$ 
est dit $n$-hyperbolique (par rapport aux 
applications harmoniques) si  pour  tout $a$ 
dans $V$, il existe un  voisinage $W$  de $a$ 
et une constante  $M>0$ tels que pour toute application harmonique 
de la boule unit\'e de ${\Bbb R}^n$ \`a valeurs dans $V$ 
et telle que $f(0)\in W$, on ait: 
$\vert f'(0)\vert  \leq M$. \\

On rappelle les notions classiques  \cite{k}
dans le cas holomorphe avec lesquelles on va faire le lien: 

\noindent 3'. Une vari\'et\'e complexe est dite 
Brody-hyperbolique s'il n'existe pas d'application 
holomorphe enti\`ere non constante \`a valeurs 
dans cette vari\'et\'e. 

\noindent Une vari\'et\'e complexe $V$  est dite 
Kobayashi-hyperbolique (ou hyperbolique) 
 si  pour  tout $a$ 
dans $V$, il existe un  voisinage $W$  de $a$ 
et une constante  $M>0$ tels que pour toute application holomorphe  
du disque unit\'e  unit\'e \`a valeurs dans $V$ 
et telle que $f(0)\in W$, on ait: 
$\vert f'(0)\vert \leq M$. 
(On pourra mettre  une m\'etrique Riemannienne 
sur $V$ pour donner un sens pr\'ecis \`a cette notion).
Comme nous n'avons pas  trouv\'e de r\'ef\'erence pr\'ecise dans la
litt\'erature, 
nous donnons en annexe finale une preuve de l'\'equivalence de 
cette d\'efinition de l'hyperbolicit\'e 
et de la d\'efinition standard.  \\

On remarquera que hyperbolique  
implique (aussi bien dans les cas 
harmoniques que complexe) Brody-hyperbolique. 
En effet si $f$ est une application 
enti\`ere non constante \`a valeurs dans $V$  on se place en $a$ 
tel que $f'(a)$ ne soit pas nul et on  consid\`ere la suite des 
$f(a+nx)$ pour montrer la non hyperbolicit\'e de $V$.\\

On fait maintenant le lien entre les notions r\'eelles et complexes. 
On va supposer  $n=2$. 

Un lien est donn\'e par la proposition suivante: 

\begin{proposition} Un  tube $\Omega$ 
est Kobayashi-hyperbolique (resp. Brody-hyperbolique)
si et seulement si $\omega$ est $2$-hyperbolique (
resp $2$-Brody hyperbolique). \end{proposition}

\emph{preuve} On utilise les relations  classiques 
entre les fonctions holomorphes et leur partie 
r\'eelle harmonique.\\

Le cas qui va nous int\'eresser est celui des 
domaines tubes de ${\Bbb C}^2$ et qui 
 d'apr\`es la proposition pr\'ec\'edente 
correspond aux valeurs $m=n=2$. Toutefois certains 
r\'esultats seront donn\'es en g\'en\'eral 
et d'autres r\'esultats pour $m$ quelquonque 
et $n=2$. 

Pour un domaine $\omega$ dans ${\Bbb R}^2$, 
on a deux cas:  Celui o\`u l'enveloppe convexe 
n'est pas tout ${\Bbb R}^2$ et celui o\`u 
cette  enveloppe est tout ${\Bbb R}^2$. 
Nous allons d'abord traiter le premier cas 
en donnant un r\'esultat complet. On remarquera 
qu'un tel domaine est contenu dans un demi-plan et  
sans nuire  \`a la  g\'en\'eralit\'e, on 
 supposera 
par la suite  que $\omega$ est contenu dans le 
demi-plan $\{(x,y)\vert y>0\}$. On dira 
alors qu'un point $a=(a_1,a_2)\in \omega$ est {\bf  born\'e}  s'il 
n'existe pas de suite $b_k$ qui tende vers $a_2$ et  tel 
que chaque  segment $[-k,k]\times \{b_k\}$ 
soit contenu dans $\omega$.

\begin{theorem} 
Soit $\omega$ un domaine de ${\Bbb R}^2$ d'enveloppe 
convexe diff\'erente de ${\Bbb R}^2$. Alors: 

1. $\omega$ est $n$-Brody hyperbolique 
si et seulement si il ne contient pas de droite affine. 

2. $\omega$ est $n$-hyperbolique si et seulement 
si tout point $a$ de $\omega$ est born\'e. 
\end{theorem}

\noindent {\bf Remarque} On voit que dans ce cas, 
$n$ ne joue aucun r\^ole.

\begin{corollary} Il existe des  tubes de 
${\Bbb C}^2$ qui sont  hyperboliques 
et dont l'enveloppe d'holomorphie n'est pas Brody- hyperbolique. 

Il existe des  tubes Brody-hyperboliques 
et non hyperboliques. \end{corollary}

Les exemples du corollaire sont faciles \`a construire
\`a partir du th\'eor\`eme.  
On va donc prouver le th\'eor\`eme. 

\emph{preuve}
\noindent 1. Il est clair que
$\Omega$ est non $n$-Brody-hyperbolique 
si $\omega$ contient une droite. 
Pour la r\'eciproque, on suppose $f=(f_1,f_2)$ 
harmonique enti\`ere \`a valeurs dans 
$\omega$. D'apr\`es 
le th\'eor\`eme de Liouville, $f_2$ 
est constante. Donc si $f$ 
est non constante, son image est une droite. 
Par cons\'equent, si $\omega$ ne contient 
pas de droite, il est $n$-Brody hyperbolique. 

\noindent 2. On suppose que $\omega$ contienne 
un point $a=(a_1,a_2)$ non born\'e. Il existe 
alors une suite $b_k$ tendant vers $a_2$ 
telle  que l'image de $]-1,1[$ par les 
applications: $t\to (kt,b_k)$ soit contenue 
dans $\omega$. De ce fait on voit imm\'ediatement que  $\omega$ est 
non $1$-hyperbolique et aussi non $n$-hyperbolique 
pour tout $n$. 

Supposons maintenant que $\omega$ soit non 
$n$-hyperbolique. Il existe alors une suite 
$f_k=(u_k,v_k)$ d'applications harmoniques
de la boule unit\'e  de ${\Bbb R}^n$ dans $\omega$ 
telles que $f_k(0)$ tende vers $a$ 
et $\vert f_k'(0)\vert $ tende vers l'infini. 
La suite $v_k$ de fonctions harmoniques 
positives est une suite normale. Quitte \`a 
extraire, on peut supposer que cette suite tend 
vers une fonction harmonique $v$ (elle ne 
peut tendre vers l'infini). La suite $u_k$ 
est alors n\'ecessairement non normale 
et on peut la renormaliser apr\`es 
extraction en une suite $U_k(x):=u_k(c_kx+d_k)$ 
qui tende uniform\'ement  sur tout compact  vers une fonction affine non constante. 
La suite des fonctions $v_k(c_kx+d_k)$ 
tend vers $a_2$ car $c_k$ et $d_k$ tendent vers $0$. 
Il est alors clair que le point $a$ 
est non born\'e. \\

On va maintenant traiter le second cas 
plus difficile \`a savoir celui o\`u 
l'enveloppe convexe de $\omega$ est tout ${\Bbb R}^2$. 
Ceci permettra d'avoir des domaines tubes hyperboliques 
dont l'enveloppe d'holomorphie est tout ${\Bbb C}^2$. 
On \'enoncera d'abord un lemme de renormalisation 
pour les applications harmoniques proche de celui 
de \cite{b} concernant les fonctions holomorphes.
Pour une application harmonique $f$, on note $\tilde{f}$ 
la somme des $\tilde{f_i}$ o\`u les 
$f_i$ sont les coordonn\'ees de $f$.  

\begin{proposition} Soit $f_k$ une suite d'applications 
harmoniques d\'efinies sur un domaine $U$  de ${\Bbb R}^n$ 
\`a valeurs dans ${\Bbb R}^m$. On suppose qu'il 
existe $p\in U$ et une suite $p_k$ qui tend vers $p$ 
telle que la suite des $\tilde{f_k}(p_k)$ tende vers 
l'infini. Alors on peut renormaliser la suite des $f_k$ 
en une suite $F_k(x):= f_k(c_kx+d_k)$ avec des $c_k>0$ 
et tendant vers $0$ et des $d_k$ tendant vers $p$ telle que:

La suite des $F_k$ (apr\`es extraction) tend uniform\'ement sur tout compact 
de ${\Bbb R}^n$ vers une application $F$ dont les 
coordonn\'ees sont des fonctions affines ou des fonctions 
identiquement $+\infty$ ou $-\infty$. De plus une 
au moins des coordonn\'ees est une fonction affine 
non constante. \end{proposition} 

\emph{preuve} On raisonne
\`a partir du lemme g\'en\'eral 
de renormalisation  comme pour les fonctions harmoniques. 
On utilise Ascoli pour montrer qu'on peut extraire des 
$F_k$ une sous-suite encore not\'ee $F_k$ 
qui tend uniform\'ement sur tout compact  vers une fonction $F$
\`a valeurs dans $[-\infty,+\infty]^m$. De plus 
la limite les $\tilde{F_k}(0)$ vaut  $1$ 
et la limite sup\'erieure des $\tilde{F_k}(a)$ 
est plus petite que $1$ pour tout $a$. On d\'eduit le r\'esultat 
en raisonnant composante par composante. On notera 
que quitte \`a extraire, on peut supposer que pour 
une des coordonn\'ees $g_k$ de $F_k$, on 
aura une limite $\tilde {g_k}(0)$ strictement positive, 
ce qui fournira une fonction limite affine non constante. \\

On donne un corollaire de la proposition 5. concernant 
les domaines  
$\omega$ de dimension deux. 

\begin{corollary} A. Un domaine $\omega$ dans ${\Bbb R}^2$  
non $n$-hyperbolique a au moins une des
deux propri\'et\'es suivante: 
 
1. Son adh\'erence contient une droite affine. 

2. pour tout $t\in{\Bbb R}$ il existe 
une suite $(x_k,y_k)$ telle que $x_k$ 
tende vers $+\infty$ et $y_k$ vers $t$ 
OU une propri\'et\'e analogue en rempla\c{c}ant 
$+\infty$ par $-\infty$ ou en permutant les coordonn\'ees.\\ 

B.Il existe des  tubes de ${\Bbb C}^2$ 
hyperboliques  dont l'enveloppe d'holomorphie est tout ${\Bbb
  C}^2$. 
\end{corollary}

\emph{preuve} A.  Si $\omega$ est non $n$-hyperbolique, 
alors on peut trouver une suite d'applications  $f_k$ 
v\'erifiant  
les hypoth\`eses de la proposition 5. \`a l'origine de 
la boule unit\'e et \`a valeurs dans $\omega$. 
On renormalise alors cette suite $f_k$ et 
on obtient une des propri\'et\'es 1 ou 2. 

B. Il est clair qu'il existe des domaines  
dans ${\Bbb R}^2$ dont l'enveloppe 
convexe est tout ${\Bbb R}^2$ 
et qui ne satisfasse \`a aucune des 
propri\'et\'es 1. et 2
(prendre par ex un voisinage 
effil\'e de trois demi-droites
convenablement choisies  partant de l'origine).  Le th\'eor\`eme de 
Bochner permet de conclure. \\

\noindent{\bf Remarque} Un domaine 
de ${\Bbb R}^2$ peut avoir la propri\'et\'e 1. 
et toutefois \^etre $n$-hyperbolique. 
Un exemple de tel domaine est donn\'e par: 
$\{(x,y)\vert 0<y<\exp(-\vert x\vert)\}$. 
Il est hyperbolique en vertu du th\'eor\`eme 6.
et clairement la droite d'\'equation 
$y=0$ est adh\'erente au domaine. 

\subsection{\bf Remarques sur les domaines non 
$n$-Brody-hyperbolique}

Nous nous pla\c{c}ons dans ${\Bbb R}^2$. 
On fait ici quelques remarques sur les images 
dans ${\Bbb R}^2$ des applications harmoniques enti\`eres. 
Le cas d\'eg\'en\'er\'e fait l'objet de la 
proposition suivante (On se place dans un cadre local). 

\begin{proposition} Soit $F$ une application 
harmonique d'un domaine $U$ de ${\Bbb R}^n$ 
\`a valeurs dans ${\Bbb R}^2$. On suppose 
$F$ de rang au plus un en chaque point 
(si $n=2$, ceci signifie que le Jacobien 
est identiquement nul). Alors $F(U)$ est 
contenue dans une droite. \end{proposition} 

\emph{preuve} Si $F$ est de rang nul 
partout, alors $F(U)$ est un point. 
On suppose donc $F$ de rang un  en un point $a$  
et donc de rang un dans un voisinage de $a$. 
Le th\'eor\`eme du rang implique alors 
qu'il existe une fonction r\'eguli\`ere 
$M$ d\'efinie sur un ouvert non vide de ${\Bbb R}$ 
et telle qu'on ait une \'equation locale: 
$u=M\circ v$ avec  $u$ et $v$  fonctions coordonn\'ees 
de $F$. Un calcul direct montre que: 
$\Delta u= M''\circ v \vert {v'}\vert ^2  +(M'\circ v) \Delta v$. 
En tenant compte de l'harmonicit\'e de $u$ et $v$ 
et du fait que $v$ est non constante, 
on en d\'eduit que $M''$ est localement nulle. 
Donc il existe deux r\'eels $c$ et $d $ 
tels que: $u=cv+d$ au voisinage de $a$. 
L'analycit\'e de $u$  et $v$ permet d'\'etendre cette
relation \`a tout $U$.\\

On est donc ammen\'e \`a \'etudier 
le cas des applications harmoniques 
enti\`eres non d\'eg\'en\'er\'ees
 dont l'image est alors d'int\`erieur non vide. 
Nous nous contenterons de quelques remarques. \\

1. Pour une fonction enti\`ere holomorphe non d\'eg\'en\'er\'ee
(c'est \`a dire non constante) le th\'eor\`eme de Picard 
nous dit que l'image est le plan ou le plan priv\'e
d'un point. Ceci est loin d'\^etre vrai pour les applications 
enti\`eres harmoniques. Dans \cite{ah} un exemple 
d'une application enti\`ere harmonique
non d\'eg\'en\'er\'ee  de ${\Bbb R}^2$ dans 
lui-m\^eme est donn\'e telles 
que les boules contenues dans l'image 
aient un rayon born\'e. Cet exemple servait \`a montrer l'inexistence 
d'une constante de Bloch. Nous donnons 
ici un exemple similaire:
On v\'erifie ais\'ement que  
l'ensemble $W:= \{(x,y)\vert 0<xy\leq 1\}\cup\{(0,0)\}$
est l'image de ${\Bbb C}$ par l'application harmonique: 
$z\to (\Re e ^z, \Re e ^{-z})$. 

Pour  cet exemple, outre l'inexistence 
de boule de rayon arbitrairement grand 
qui y soit contenue, on notera qu'il ne contient pas de 
droite mais que deux droites y sont adh\'erentes. 
A partir de l\`a, il est  facile de construire des 
exemples de domaines non   $2$-Brody hyperbolique
qui ne  contiennent  pas de droite 
et qui sont de ce fait  $1$-Brody hyperbolique.

La proposition suivante caract\'erise, \`a transformation 
affine pr\`es, 
les  applications harmoniques enti\`eres qui sont 
holomorphes.  

\begin{proposition} Une application harmonique $H$
enti\`ere de ${\Bbb R}^2$
dans lui-m\^eme  dont  le Jacobien est positif 
ou nul en tout point  
 est \`a transformation affine pr\`es
(sur l'ensemble image)  une 
fonction holomorphe  enti\`ere.\end{proposition} 

On a alors le corollaire: 

\begin{corollary} 
 Une application harmonique $H$
enti\`ere de ${\Bbb R}^2$
dans lui-m\^eme dont le 
Jacobien ne s'annule pas,  
 est \`a transformation affine pr\`es
(sur l'ensemble image)  une 
fonction holomorphe  enti\`ere.\end{corollary}

\emph{preuve de la proposition} Soient $u$ et $v$ les coordonn\'ees de 
$H$. On pourra \'ecrire: 
$u=\Re f$ et $v=\Re g$ avec $f$ et $g$ des fonctions 
enti\`eres. Un calcul direct montre 
que le Jacobien de $F$ est donn\'e par 
$\Im (f'\overline{g'})$. En \'eliminant 
le cas $g$ constant qui est trivial, 
la positivit\'e de ce jacobien 
signifie que $\Im (\frac{f'}{g'})$ 
est positif  l\`a o\`u la fonction m\'eromorphe 
$\frac{f'}{g'}$ est d\'efinie. Par une forme 
du th\'eor\`eme de Liouville, ceci implique que 
$\frac{f'}{g'}$ est constant. Des consid\'erations 
\'el\'ementaires permettent de conclure. \\

\section{\bf Applications  \`a valeurs dans un tore
ou un groupe de Lie}

Dans ce chapitre, on consid\`ere 
des applications solutions d'une 
edp elliptique \`a valeurs dans un tore. 
Concernant le Laplacien, les r\'esultats de 
renormalisation sont d'une certaine 
fa\c{c}on plus faibles 
que ceux pr\'ecedemment obtenus. 
Pour ce qui concerne les op\'erateurs 
elliptiques, on renvoie au livre 
de H\"{o}rmander \cite{ho}. 

Dans la suite, on consid\`ere  
 un polyn\^ome homog\`ene r\'eel  $P$ 
sur ${\Bbb R}^m$ auquel on associe  naturellement 
un op\'erateur diff\'erentiel $P(\partial)$. 
On va supposer que $P$ ne s'annule qu'en 
$0$ sur ${\Bbb R}^m$. L'op\'erateur 
$P(\partial)$ est alors elliptique. Rappelons 
une propri\'et\'e essentielle: 

Soit $U$ un domaine de ${\Bbb R}^m$. 
Toute fonction $f$ suffisemment r\'eguli\`ere  sur $U$ 
et solution de $P(\partial)f=0$ 
est analytique r\'eelle.
On appellera $P$-fonction une telle fonction.
On a la propri\'et\'e suivante: Si  une suite de $P$-fonctions 
d\'efinies sur $U$ converge uniform\'ement sur tout compact 
vers une fonction $f$, alors $f$ est aussi une $P$-fonction. 
De plus les d\'eriv\'ees convergent aussi uniform\'ement sur tout
compact. 

Introduisons les d\'efinitions suivantes: 

1. Une application de $U$ dans ${\Bbb R}^n$ 
est une $P$-application si ses coordonn\'ees 
sont des $P$-fonctions. 

2.  Une application $F$  de $U$ dans $T_n:={\Bbb R}^n/{{\Bbb Z}^n}$ 
est une $P$-application si elle s'\'ecrit sous la forme: 
$p\circ G$  avec $p$ projection canonique 
de ${\Bbb R}^n$ sur $T_n$ et $G$ une $P$-application. 
Notons que $G$ est d\'efinie \`a une constante additive 
pr\`es. De ce fait $G'$ qu'on notera $F'$ 
dans la suite, est bien d\'efinie. C'est une 
$P$-application. 

On a le th\'eor\`eme suivant de renormalisation. 

\begin{theorem} Soit $f_k$ une suite 
de $P$-application de $U$ \`a valeurs 
dans $T_n$. On suppose que pour 
un point $p$ de $U$, il existe une suite 
$p_k$ tendant vers $p$ et telle que 
$\vert f_k'(p_k)\vert$ tende vers l'infini. 
Alors on a une renormalisation d'une sous-suite 
des $f_k$ en une fonction affine non constante
(notons que cette derni\`ere notion a un sens).
La renormalisation signifie qu'il existe une suite 
de nombres positifs  
$a_k$ tendant vers $0$  et une suite de vecteurs $b_k$ 
tendant vers $p$ 
et telles qu'une suite extraite des 
$f_k(a_kx+b_k)$ 
converge sur tout compact 
de ${\Bbb R}^m$  vers une fonction affine. \end{theorem} 

Le th\'eor\`eme implique un corollaire 
analogue \`a celui des fonctions 
harmoniques concernant les $P$-applications enti\`eres 
. Le corollaire 
qui se d\'eduit imm\'ediatement 
du th\'eor\`eme,  fait le lien avec la renormalisation \`a 
constante pr\`es.

\begin{corollary}        
 Soit $f_k$ une suite 
de $P$-application de $U$ \`a valeurs 
dans ${\Bbb R}^n$. On suppose les 
m\^emes hypoth\`eses que dans le th\'eor\`eme 
pour un point $p$. Alors il existe une suite 
de vecteurs $c_k$ de ${\Bbb R}^n$, 
telle qu'on puisse renormaliser 
une sous suite  des $f_k+c_k$ en une 
application affine non constante. \end{corollary}

\noindent {\bf Remarques} Avant de d\'emontrer 
le th\'eor\`eme, faisons quelques remarques 

1. Le th\'eor\`eme similaire de M. Green  concernant 
les applications holomorphes \`a valeurs dans 
un tore complexe est bien connu. En fait, 
il est facile de retrouver ce r\'esultat 
\`a partir de n\^otre th\'eor\`eme 
en consid\'erant une fonction holomorphe comme un couple de 
fonctions harmonique v\'erifiant les relations de 
Cauchy-Riemann.  

2. La notion de renormalisation \`a constante 
pr\`es concernant le cas holomorphe \`a une variable 
est \`a  la base du formalisme de Minda
(voir le livre \cite{s} et l'article \cite{m}).  
Cette renormalisation  est essentiellement 
\'equivalente \` a  la renormalisation des 
fonctions holomorphes dans un tore complexe.

Il nous parait important 
de remarquer  que dans le cas des fonctions harmoniques  
 les r\'esultats 
de renormalisation obtenus par 
cette m\'ethode sont plus faibles que ceux d\'emontr\'es 
pr\'ec\'edemment.  
En g\'en\'eral une $P$-fonction enti\`ere  \`a 
valeurs dans ${\Bbb R}$ n'est pas renormalisable 
en une fonction affine non constante 
comme le montre l'exemple simple des $P$-fonctions  
que sont les  polyn\^omes \`a une variable. \\  

\emph{Preuve du th\'eor\`eme}
Pour une $P$-application $F$, on consid\`ere 
$\vert F'\vert $ qui joue le m\^eme r\^ole 
que $\tilde{f}$. 

 Le lemme 
 de renormalisation  
nous dit qu'il existe deux suites comme dans le th\'eor\`eme 
telle que  pour la suite des 
$s_k(x):=f_k(a_kx+b_k)$ , on a: 
la suite des  $s'_k(x)$
tend uniform\'ement sur tout compact vers une limite 
plus petite que $1$ en norme  et telle que les $s'_k(0)$ 
valent $1$ en norme. Un argument de normalit\'e 
valable dans le cas elliptique 
permet d'affirmer qu'une sous-suite de 
$s'_k$ a une limite. Cette limite $t$ est une 
$P$-application {\bf born\'ee}. Le 
th\'eor\`eme de Liouville (voir lemme suivant) 
nous dit que $t$ est constante.  Pour conclure, on 
peut utiliser 
le th\'eor\`eme classique d'int\'egration 
des limites de suites de fonctions jointe \`a 
la compacit\'e de $T_n$. \\

\begin{lemma} Une $P$-fonction $f$  born\'ee est 
constante. \end{lemma}

\emph{preuve} Par l'ellipticit\'e 
de $P(\partial)$, la transform\'ee de 
Fourier de la distribution temp\'er\'ee 
$f$ est \`a support en l'origine, 
donc une combinaison d'une masse de Dirac et ses d\'eriv\'ees 
en l'origine. Par cons\'equent $f$ est un polyn\^ome 
et comme $f$ est born\'ee, c'est une constante.

\subsection{le cas des groupes de 
Lie complexes}

Les id\'ees pr\'ec\'edentes 
peuvent aussi s'appliquer aux groupes  de Lie complexes.
Nous ne tra\^iterons pas ce sujet dans toute sa g\'en\'eralit\'e.  
On fait d'abord quelques rappels: 

Soit $G$ un groupe de Lie complexe 
(le lecteur peu famillier avec cette notion
pourra se placer dans la situation 
$G=Gl(n, {\Bbb C})$). 
Pour une application  holomorphe $f$ 
d'un domaine  $U$ de ${\Bbb C}$ \`a 
valeurs dans $G$, on pose: 
$Df(z):= f(z)^{-1}f'(z)$. 
L'application $Df$ est holomorphe 
\`a valeurs dans l'espace tangent 
en l'idendit\'e (alg\`ebre de Lie), qui est isomorphe 
\`a  un ${\Bbb C}^n$. Si 
$Df(z)$ est constante, alors 
$f(z)$ est de la forme $g\exp (zX)$ pour un certain 
$g\in G$ et un certain  $X$ dans l'alg\`ebre de Lie 
de $G$. Ces applications g\'en\'eralisent les
applications affines.  

On a le th\'eor\`eme suivant: 

\begin{theorem} Soit $f_k$ une suite d'applications 
holomorphes de $U$ dans $G$. On suppose 
qu'il existe $p\in U$ tel que la norme 
de $Df_k(p)$ tende vers l'infini. Alors 
il existe une renormalisation des 
$f_k$ \`a constante pr\`es  en une fonction enti\`ere de 
la forme $g\exp (zX)$. avec $X$ non nul.  De mani\`ere explicite: 
Il existe une suite $g_k$ dans $G$, 
une suite  de nombres positifs $a_k$ tendant vers $0$ 
et une suite $b_k$ d'\'el\'ements de $U$ tendant vers $p$ 
tels que la suite de fonctions $g_kf_k(a_kz+b_k)$ 
tende uniform\'ement vers $g\exp (zX)$ sur 
tout compact de ${\Bbb C}$. \end{theorem} 

\emph{preuve} La preuve utilise les id\'ees 
pr\'ec\'edentes.  On applique  le lemme g\'en\'eral 
de renormalisation aux $\vert Df_k\vert$. 
Ceci permet d'obtenir une suite renormalis\'ee 
de fonctions $u_k(z):=f_k(a_kz+b_k)$.
On remplace ensuite la suite 
des $u_k$ par une suite $U_k(z):=g_k u_k(z)$ 
avec les $g_k$ choisis 
de telle sorte que $U_k(p)= id$. On 
a $DU_k=Du_k$. La suite des $U_k$ 
est alors normale et elle tend 
apr\`es extraction vers une 
application holomorphe  enti\`ere non constante $F$
\`a valeurs dans $G$. La fonction $F$ 
est de Brody, c'est \`a dire qu'elle v\'erifie 
$\vert DF(z)\vert \leq 1$. Le th\'eor\`eme de 
Liouville classique nous dit alors que $DF$ 
est constante. Le th\'eor\`eme se d\'eduit 
des remarques de d\'ebut  sur les groupes de lie. \\

{\bf Remarques}

1.Comme cas particulier du th\'eor\`eme, 
on a la renormalisation \`a constante 
pr\`es d'applications holomorphes 
enti\`eres non constantes. 

2. On a le m\^eme type de 
r\'esultats si on consid\`ere des 
applications holomorphes \`a valeurs 
dans un quotient de $G$ 
par un groupe discret cocompact. 
Toutefois
 la situation 
${\Bbb C}^{\star}={\Bbb C}/{\Bbb Z}$ 
\'etudi\'ee par Berteloot et Duval,
est en relation avec la renormalisation 
\`a constante {\bf r\'eelle} pr\`es.
Dans le cas d'un groupe de Lie complexe, 
ceci signifie 
une renormalisation \`a une constante 
pr\`es variant dans une forme r\'eelle. 
Une telle renormalisation 
 pourrait \^etre utilis\'ee 
pour l'\'etude de tubes g\'en\'eralis\'es,i.e. des domaines
invariants par l'action de la  forme r\'eelle ( Voir 
\cite{c} ). 

3. On peut g\'en\'eraliser la notion de $P$- 
application pour un groupe de Lie r\'eel $G$ 
de la mani\`ere suivante: Une application 
$f$ d'un domaine $U$ de ${\Bbb R}^n$
\`a valeurs dans $G$ est une $P$-application 
si $P(Df)=0$. On peut alors d\'evelopper 
une th\'eorie similaire \`a celle du cas Euclidien. 
Nous ignorons toutefois la port\'ee
que peut avoir cette notion.

\section{Annexe sur l'hyperbolicit\'e}

On note $D_r$ le disque  ouvert centr\'e 
\`a l'origine de rayon $r>0$. 

\begin{proposition} 
Soit $X$ une vari\'et\'e complexe qu'on munit d'une 
m\'etrique Riemannienne. 
Alors $X$ est hyperbolique au sens de Kobayashi 
si et seulement si on a la propri\'et\'e (P)  suivante: 

Pour tout point $a$ de $X$, il existe un voisinage $V$ 
et une constante $M>0$ tels que pour toute application 
holomorphe $f:D_1\to X$ telle que $f(0)\in V$, on ait:  
$\vert f'(0)\vert \leq M$. \end{proposition} 

\noindent\emph{preuve}
On suppose d'abord $X$ hyperbolique. Pour $a$ dans 
$X$ on choisit une boule ferm\'ee $B$
pour la m\'etrique de Kobayashi centr\'ee en $a$ 
et de rayon $2r>0$ assez petit pour que cette boule soit compacte. 
(Ceci est possible car la m\'etrique de Kobayashi induit la topologie 
\cite{k}). On peut alors choisir  pour $V$ la boule 
ferm\'ee de centre $a$ et de rayon $r$ pour la 
m\'etrique de Kobayashi. En effet d'apr\`es la 
propri\'et\'e de contraction de la m\'etrique  
 de Kobayashi, il existe $r'>0$ 
tel que pour toute application holomorphe 
$f:D_1\to X$ v\'erifiant $f(0)\in V$, on 
ait: $f(D_{r'})\subset B$. Comme $B$ 
est compact, un argument de normalit\'e 
permet de montrer l'existence de $M$. 

La r\'eciproque est plus difficile. On suppose 
la propri\'et\'e (P) v\'erifi\'ee pour $X$. 
En chaque point $p$, on d\'efinit 
l'indicatrice $K_p$  de Kobayashi comme 
\'etant le sous-ensemble de l'espace 
tangent $T_p$ form\'e des \'el\'ements $rf'(0)$ 
avec $0<r<1$ avec  $f$ application holomorphe 
de $D_1$ dans $X$  v\'erifiant $f(0)=p$.
A cette indicatrice, on associe de 
mani\`ere classique une jauge $j_p$ 
d\'efinie sur $T_p$ par: 
$j_p(v)= \inf \{t>0\vert v/t \in K_p\}$.
La propri\'et\'e (P) signifie 
alors simplement que pour tout point $a$, il existe 
un voisinage $V$ de $a$ et $c>0$ tel 
que pour tout $p\in V$ et $x\in T_p$, on a: 
$j_p(x)\geq c\vert x\vert$. Le point essentiel 
qui permet de conclure est maintenant 
que la m\'etrique de Kobayashi 
est la m\'etrique int\'egr\'ee par rapport 
aux $j_p$ (voir \cite{k} pour une preuve 
de ce th\'eor\`eme assez difficile).  \\

\end{document}